# GENRE ET EVALUATION EN MATHEMATIQUES : ETUDE COMPARATIVE DE GESTION D'EPISODES EVALUATIFS

## BRISMONTIER[1*] Chloé

**Résumé –** L'article s'intéresse aux différences de performance en mathématiques entre filles et garçons visibles dès les quatre premiers mois de scolarisation obligatoire dans le système éducatif français. L'influence des stéréotypes de genre dans les pratiques évaluatives des professeur·es et la menace du stéréotype de genre sur la performance des élèves y sont questionnées. Pour obtenir des éléments de réponse, une étude comparative de gestion d'épisodes évaluatifs est proposée à partir d'outils théoriques.

**Mots-clés** : épisodes évaluatifs, stéréotypes de genre, régulation, performance, affect

**Abstract –** The article focuses on the differences in mathematics performance between girls and boys visible from the first four months of compulsory schooling in the French education system. The influence of gender stereotypes in the evaluation practices of teachers and the threat of the gender stereotype on student performance are questioned. To obtain answers, a comparative study of management of evaluative episodes is proposed based on different theoretical tools.

## I. INTRODUCTION

En France, bien que les évaluations nationales aient pour ambition d'offrir une approche édumétrique et psychométrique précise, elles n'apportent aucune explication aux causes de l'écart de performance en mathématiques entre les filles et les garçons, un écart souvent défavorable aux filles. Cette absence d'éclaircissements fragilise la lutte scientifique contre les hypothèses controversées sur les différences intellectuelles entre les genres – i.e. sexe biologique – (C. Perronet et al., 2024). Plusieurs disciplines, telles que les sciences de l'éducation, la psychologie et la sociologie, ont proposé diverses explications. Parmi elles, les pratiques évaluatives des professeur·es sont aujourd'hui questionnées. Elles sont soupçonnées d'être influencées par des stéréotypes de genre, modifiant l'affect mathématique des élèves, au sens de la théorisation de V. A. Debellis et G. A. Goldin (2006), notamment les croyances, les valeurs, les émotions et les représentations liées aux mathématiques, et serait préjudiciable quant à la performance des filles. Cela pourrait se traduire par un contrat didactique en évaluation différencié selon le genre de l'élève. Nous entendons par contrat didactique en évaluation un contrat didactique spécifique à un épisode évaluatif proposé par un·e professeur·e à ses élèves au cours de l'étude d'un savoir (Sayac, 2017). Cette notion de contrat didactique différencié constitue l'objet central de ma thèse, où j'en examine les caractéristiques ainsi que ses liens avec la performance et l'affect mathématique des élèves. Dans cet article, je propose d'abord une analyse graphique des résultats des élèves aux évaluations nationales de la direction de l'évaluation, de la prospective et de la performance (DEPP), mettant en lumière l'écart de performance entre filles et garçons en mathématiques. Ensuite, je passe en revue les principaux travaux scientifiques portant sur l'interprétation de cet écart, avant de présenter une partie du cadre de mon étude sur l'influence du genre dans les pratiques évaluatives des professeur·es, les relations entre celles-ci, les performances scolaires et l'affect des élèves envers les mathématiques.

---

[1*] Université Paris Cité – France – chloe.brismontier@etu.u-paris.fr

## II.    LES ÉVALUATIONS NATIONALES DE LA DEPP

Depuis l'année scolaire 2018-2019, la DEPP en collaboration avec la direction générale de l'enseignement scolaire (DGESCO) et l'inspection générale de l'éducation, du sport et de la recherche (IGESR), a instauré un dispositif d'évaluation nationale dénommé « Évaluations Repères ». Ce dispositif a pour objectif d'évaluer les niveaux de maîtrise en français et en mathématiques des élèves dans les écoles françaises, tout en permettant aux professeur·es d'ajuster l'organisation pédagogique de leurs classes. Les évaluations de mathématiques destinées aux élèves de cours préparatoire (CP), cours élémentaire 1ère année (CE1), cours moyen 1ère année (CM1) et 6ème année (6ème) - respectivement grade 1, grade 2, grade 4 et grade 6 pour les niveaux scolaires anglophones - constituent l'une des principales sources de données pour mon étude.

Dans le cadre de ma thèse en didactique des mathématiques, je m'intéresse spécifiquement à l'écart de performance entre les filles et les garçons dans les écoles élémentaires (du grade 1 au grade 5) et les classes de 6e, soit pour les élèves âgés de 6 à 12 ans. Je propose de commencer par une analyse des données sur les résultats des « Évaluations Repères » par niveaux scolaires selon les types de tâches mathématiques de 2018 à 2023. Puis, de poursuivre par une analyse graphique des résultats par domaine mathématique.

*1.   Résultats aux évaluations nationales par niveaux scolaires selon le type de tâche mathématique de 2018 à 2023*

La première vague du dispositif d'évaluation nationale de la DEPP s'est déroulée en septembre 2018. Au départ, les « Évaluations Repères » concernaient uniquement les élèves inscrit·es en CP, CE1 et 6e. Par la suite, le dispositif a été étendu aux classes de CM1 à partir de l'année scolaire 2023-2024, et aux classes de CE2 et CM2 à partir de l'année scolaire 2024-2025. À l'exception des élèves de CP, chaque élève est évalué·e lors d'une passation unique en début d'année scolaire. Pour les élèves de CP, deux passations sont organisées : une en début d'année et une autre en milieu d'année scolaire. À la suite de chaque passation des « Évaluations Repères », la DEPP publie des notes d'information présentant les résultats issus du traitement statistique des réponses des élèves. Ces notes d'information, disponibles sur le site de la DEPP[2], sont accompagnées d'un fichier contenant les données à partir desquelles les tableaux et les graphiques ont été élaborés. Grâce aux données mises en ligne par la DEPP, il est possible d'avoir une lecture de l'écart de performance entre filles et garçons par niveau scolaire selon le type de tâche mathématique évalué pour certaines évaluations nationales.

D'après les données de la DEPP pour le CP, l'écart de performance entre les filles et les garçons est en faveur des filles pour presque tous les types de tâche mathématique évalués lors des passations en début de CP en septembre 2018, 2019, 2020 et 2023. De 2018 à 2022, les filles inscrites en CP sont moins performantes que les garçons en ce qui concerne la comparaison des nombres et le placement d'un nombre sur une ligne graduée. En revanche, pour la passation de début CP en 2023, les filles sont moins performantes que les garçons pour la reproduction d'un assemblage et la résolution de problèmes.

---

[2] https://www.education.gouv.fr/les-evaluations-nationales-exhaustives-307627

Un phénomène particulier apparaît à la mi-CP. Les résultats de chaque passation des « Évaluations Repères » depuis 2020 montrent que les filles à la mi-CP sont moins performantes que les garçons pour l'ensemble des types de tâche évalués en mathématiques.

Pour les élèves inscrits en CE1, les écarts de performance entre filles et garçons, en défaveur des filles, sont beaucoup plus importants et nombreux qu'en début de CP. Le phénomène observé à la mi-CP se confirme en début de CE1. Jusqu'en 2022, les filles sont plus performantes que les garçons pour la reproduction d'un assemblage et le calcul mental (et seulement pour cela). En 2023, les filles sont plus performantes que les garçons pour l'addition et le calcul mental.

D'après les données fournies par la DEPP pour les élèves inscrits en CM1 en septembre 2023, la proportion de filles atteignant une maîtrise satisfaisante est plus élevée que celle des garçons uniquement pour la pose de calculs et le calcul de nombres. Cet écart de performance en faveur des filles est minime par rapport aux écarts de performance en faveur des garçons dans les autres types de tâche évalués en mathématiques.

Pour les évaluations nationales de début de 6e, la DEPP n'a pas publié les données de l'année 2018, et les données disponibles pour les années suivantes ne répartissent pas les performances selon le genre des élèves à l'exception de la passation de septembre 2023.

2. *Résultats aux évaluations nationales par domaine mathématique selon le niveau scolaire de 2018 à 2023*

Pour mener une analyse historique des résultats des différentes cohortes d'élèves aux évaluations nationales de CP, CE1, CM1 et 6e, je propose une lecture graphique de l'évolution de la proportion de filles et de garçons atteignant une maîtrise satisfaisante dans chaque domaine évalué en mathématiques.

Je précise que les types de tâche mathématique regroupés sous « Nombres et Calculs », en plus de ceux évalués en 6e sous cette appellation, sont les suivants :

- Comparer des nombres entiers (CP),

- Lire des nombres entiers (CP et CE1),

- Écrire des nombres entiers ; Placer un nombre sur une ligne graduée (CP, CE1 et CM1),

- Additionner ; Calculer mentalement ; Soustraire ; Représenter des nombres entiers (CE1),

- Poser et calculer (CM1).

Pour « Espace et Géométrie », les types de tâche mathématique évalués en 6e sous cette appellation, ainsi que ceux des autres niveaux scolaires, sont :

- Quantifier des collections (CP),

- Reproduire un assemblage (CP et CE1).

En ce qui concerne « Résolution de problèmes », seule la résolution de problèmes évaluée en début de CP, CE1, CM1 et 6e est prise en compte. Pour « Automatismes », seuls les types de tâche mathématique évalués en 6e et la mémorisation des faits numériques et des procédures en CM1 pour l'année 2023 sont disponibles. En raison du manque de données, il n'est pas possible de fournir un graphique historique pour la proportion d'élèves filles et garçons en maîtrise satisfaisante dans le domaine « Automatismes ». Ce domaine est donc exclu de mon analyse. Enfin, le domaine « Espace et Géométrie » n'est pas évalué lors des évaluations nationales en CM1.

Les courbes des graphiques pour les domaines « Nombres et Calculs » (graphique 1) et « Résolution de problèmes » (graphique 3) présentent des similitudes. Chaque année où les données étaient disponibles, on observe une meilleure performance des élèves filles et garçons en début de CP par rapport au CE1. En 2023, une diminution de la performance est notée du CP jusqu'à la 6e. La courbe représentant la moyenne des filles en maîtrise satisfaisante dans « Nombres et Calculs » (graphique 1) et dans « Résolution de problèmes » (graphique 3) est inférieure à celle des garçons, mais suit une trajectoire similaire. Un écart moyen de -3,3 % est observé entre la courbe des élèves filles et celle des élèves garçons en « Nombres et Calculs », tandis que l'écart moyen est de -3,03 % en « Résolution de problèmes ».

Les courbes du graphique pour le domaine « Espace et Géométrie » (graphique 2) présentent une trajectoire différente de celles des deux domaines précédents. La proportion d'élèves garçons et d'élèves filles en maîtrise satisfaisante dans « Espace et Géométrie » varie très peu au fil des années et entre les niveaux de CP et CE1, jusqu'en 2022. De plus, la proportion de filles atteignant une maîtrise satisfaisante est généralement supérieure à celle des garçons, sauf pour les élèves filles inscrites en CE1 en 2023. Les résultats de la passation de 2023 montrent une forte baisse de performance du CP au CE1, et du CE1 jusqu'à la 6e.

*3. Conclusion*

La diminution des performances des élèves en mathématiques du CP jusqu'en 6e est observable, quel que soit le genre de l'élève et le domaine mathématique évalué pour les élèves en maîtrise satisfaisante. Cependant, on note une meilleure performance des filles en mathématiques en début de CP, contrairement aux résultats en CE1, CM1 et 6e, à l'exception du domaine « Espace et Géométrie ». D'un point de vue didactique, cette observation soulève des questions quant à la validité des évaluations nationales (Grapin et Mounier, 2020). En effet, la DEPP précise dans sa documentation que les évaluations en CP, CE1 et CM1 ne mesurent pas les compétences mathématiques des élèves et qu'elles n'ont pas pour vocation à mesurer tout ce qui a été appris les années précédentes, mais leurs capacités cognitives fondamentales, contrairement à celles destinées aux élèves de 6e. À la lumière de cette précision, on peut remettre en question le niveau de maîtrise des élèves de primaire en mathématiques décrit par les résultats des « Évaluations Repères ». Toutefois, ces résultats révèlent une différence non négligeable entre la proportion de filles et celle de garçons classés en maîtrise satisfaisante par la DEPP pour un grand nombre de types de tâche évalués en mathématiques.

Pour expliquer la différence de performance entre les filles et les garçons aux « Évaluations Repères », plusieurs pistes ont été explorées par divers champs disciplinaires. Je propose une synthèse, non exhaustive, des notes d'information, articles, ouvrages et enquêtes sur les caractéristiques de cet écart.

## III. ÉTAT DE LA LITTERATURE NON EXHAUSTIF SUR LES CARACTERISTIQUES DE LA DISPARITÉ DES PERFORMANCES FILLES-GARÇONS EN FRANCE

Récemment, des économistes de l'Institut des Politiques Publiques (IPP) ont publié une note intitulée « Le décrochage des filles en mathématiques » (2024). Cependant, comme le montrent les données fournies par la DEPP, la baisse de performance en mathématiques ne concerne pas uniquement les filles, mais également les garçons. Néanmoins, ces recherches de l'IPP ont permis de mieux caractériser l'écart de performance en mathématiques entre filles et garçons, en défaveur des filles. Tout d'abord, l'évolution de cet écart toucherait principalement les filles les plus performantes en début de CP et se manifesterait dans toutes les catégories sociales et configurations familiales, ainsi que sur l'ensemble du territoire. Ensuite, quelle que soit leur performance initiale, les filles progresseraient moins que les garçons de même niveau au départ. Enfin, l'écart serait moins important dans les classes où la majorité des élèves sont des filles, lorsqu'une femme est l'enseignante plutôt qu'un homme, et lorsque l'école est située dans une zone de Réseau d'Éducation Prioritaire Plus (REP+).

À travers l'enquête ELFE (Étude Longitudinale Française depuis l'Enfance), l'Institut National d'Études Démographiques (Ined) avait dès 2022 relevé des caractéristiques concernant l'écart de performance en mathématiques entre filles et garçons. Le communiqué de presse rédigé par Fischer et Thierry (2022) précise que cet écart en défaveur des filles serait détectable dès l'âge de 5 ou 6 ans, indépendamment de la classe. Les auteurs suggèrent que l'influence du genre ne serait pas d'origine innée, mais demeure inexplicable. D'autres chercheurs et chercheuses se sont intéressé·es à l'hypothèse de l'influence du genre sur la performance des élèves dans le système éducatif français. Plusieurs chercheuses, notamment en sciences de l'éducation, en sociologie et en didactique des mathématiques, remettent en question les pratiques évaluatives des professeur·es (Mathieu-Benmerah, 2023 ; Perronnet et al., 2024 ; Sayac, 2024), affirmant qu'elles reproduisent des stéréotypes de genre. Il est probable que les retours et appréciations des professeur·es soient différenciés selon le genre de l'élève, au détriment des filles. En conséquence, le sentiment de compétence et d'efficacité personnelle des filles pourrait être affaibli (Imberdis et al., 2021). En plus de la connaissance de soi et de l'estime de soi (Martinot, 2001 ; Bishop, 2014 ; Darragh et Radovic, 2018), la performance et les croyances motivationnelles joueraient également un rôle dans la réussite scolaire en mathématiques des élèves. « Les élèves intègrent certains schémas de groupes associés aux disciplines » (Jamain, 2019, p.151), étroitement liés à leur sentiment de compétence et à l'intérêt porté à chaque discipline (Imberdis et al., 2021). N. Mella (2021) confirme ce point en affirmant que la réussite scolaire est liée aux compétences socio-émotionnelles et aux mesures d'adaptation scolaire. Elle pose l'hypothèse, non vérifiée, que les relations entre l'autorégulation scolaire et les compétences socio-émotionnelles dépendent de la capacité des élèves à déployer des stratégies d'apprentissage efficaces. De plus, elle précise que la perception de l'enjeu évaluatif diffère selon la discipline et que cette perception en mathématiques pourrait être défavorable à la réussite scolaire des filles (Mella, 2021).

Une étude qualitative portant sur l'écart de performance en mathématiques entre les filles et les garçons, au désavantage des filles, s'avère nécessaire pour mieux comprendre les liens de dépendance entre la réussite scolaire, l'affect mathématique des élèves et les régulations de chaque individu. Afin d'examiner ces relations, j'ai choisi de travailler principalement dans le cadre du contrat didactique en évaluation (Sayac, 2017). Ce cadre propose un changement de paradigme qui s'ancre dans la réalité des classes en prenant en compte à la fois les facteurs conjoncturels et personnels, dans le but de décrire et d'expliquer les phénomènes relatifs aux interactions entre l'enseignement, les apprentissages et leur évaluation. Je formule l'hypothèse qu'un contrat didactique en évaluation différencié existe en fonction du genre de l'élève.

IV.    DÉFINITION DES OUTILS THEORIQUES DE MA RECHERCHE

Pour analyser l'influence du genre sur le contrat didactique en évaluation, cette recherche repose sur une étude comparative de gestion d'épisodes évaluatifs informels observés lors de séances d'apprentissage en « Résolution de problème » et en « Espace et Géométrie ». Le choix de ces deux domaines mathématiques s'appuie sur la disparité des performances entre filles et garçons, illustrée par les graphiques 2 et 3. Afin de recueillir des données en quantité suffisante, l'étude est menée dans sept écoles élémentaires et trois collèges figurant dans le territoire expérimental du projet Maryam Mirzakhani[3] du rectorat de l'Académie de Lille. L'objectif est d'examiner comment les régulations interactives et les jugements professionnels et didactiques en évaluation (JPDE) s'intègrent aux stratégies d'apprentissage des élèves et varient selon plusieurs facteurs : la tâche mathématique proposée, son degré de complexité, les pratiques évaluatives du professeur·e et le genre de l'élève. Pour ce faire, quatre niveaux d'analyse sont envisagés : l'épisode évaluatif à travers le prisme du contrat didactique en évaluation, l'activité mathématique des élèves à partir de la théorie de l'activité, les stratégies de régulation des élèves selon la conceptualisation de L. Allal (2007) et les régulations interactives entre le ou la professeur·e et ses élèves. Cet article se concentre sur ce quatrième niveau d'analyse en détaillant les types de régulations interactives observées dans les épisodes évaluatifs, en lien avec les tâches mathématiques proposées.

*1.    Catégorisation des régulations interactives*

Les épisodes évaluatifs informels étudiés correspondent à des séances d'apprentissage où sont observées les interactions quotidiennes entre le ou la professeur·e et ses élèves. Je considère que ces interactions reposent sur la capacité ou l'action du professeur·e à formuler un jugement :

- sur les connaissances mathématiques, les procédures et le résultat des élèves, en référence à la tâche mathématique assignée ou aux savoirs prescrits ou enseignés ;
- à caractère injonctif, en lien avec les éléments de contexte de l'épisode évaluatif ;
- sur l'autorégulation de l'élève ou sur l'élève lui-même.

Dans le premier cas, les jugements du professeur·e apportent une information qui agit sur les processus ou les dimensions de l'activité mathématique de l'élève, respectivement lors d'une activité en résolution de problème comme modélisé par N. Grapin et É. Mounier (2024) ou en géométrie comme définis par T. Barrier, C. Hache, A.-C. Mathé et S. Montigny (2013). Dans le deuxième cas, les jugements du professeur·e apportent une information quant à ses JPDE. Pour finir, dans le troisième cas, les jugements du professeur·e apportent une information qui agit sur l'affect mathématique de l'élève et sa capacité à mettre en place des stratégies d'auto-régulation (Debellis et Goldin, 2006). De ma description des interactions quotidiennes pouvant déboucher sur des régulations interactives à partir de la capacité ou l'action du professeur·e à porter des jugements s'en dégage trois registres de régulations interactives : un registre didactique, un registre injonctif et un registre affectif. À partir de cette catégorisation des régulations interactives, des critères d'analyse de celles-ci seront développés pour permettre l'étude de leurs régularités et de leurs variabilités, afin d'identifier les facteurs conjoncturels et personnels à l'œuvre dans la gestion d'épisodes évaluatifs en résolution de problème et en géométrie.

---

[3] https://filles-maths-nsi-projet-maryam-mirzakhani.site.ac-lille.fr

## 2. Analyse des régularités et variabilités des régulations interactives

En s'appuyant sur les travaux de M. Kiwan-Zacka et E. Roditi (2018), sur l'analyse des déroulements, qui associent un couple information-action pour chaque « régulation didactique » entre le ou la professeur·e et l'élève ou les élèves, je vais construire mes critères d'analyse des régulations interactives en déterminant quelles informations et quelles actions peuvent avoir lieu dans chaque registre.

Pour le registre didactique, les couples information-actions possibles sont délimités à partir de la définition de « régulations didactiques » de M. Kiwan-Zacka et E. Roditi (2018).

> « Une « régulation didactique » est caractérisée par la production d'un ou plusieurs élèves lors de la réalisation d'une tâche mathématique, de son interprétation par l'enseignant, et de l'intervention de ce dernier qui peut être un aide à la réalisation de la tâche ou une modification – simplification – de la tâche à réaliser » (Kiwan-Zacka et Roditi, 2018).

Afin de distinguer une action comme aide ou tentative d'aide à la réalisation d'une tâche mathématique d'une action de modification de la tâche à réaliser pour les deux domaines mathématiques observés. Je m'appuie sur la catégorisation des aides apportées dans une construction instrumentée de E. Petitfour (2015) et sur la construction d'une catégorisation des aides à la résolution de problème à partir d'articles d'écrivant les actions des professeur·e pour venir en aide à leurs élèves (Croguennec, 2023; Veyrunes et al., 2005). Les couples information-actions possibles d'après ce travail théorique concernent les combinaisons entre les informations portant sur les connaissances, la procédure ou le résultat de l'élève et les actions amenant à une modification de la tâche mathématique ou à une aide ou une tentative d'aide à la résolution de la tâche mathématique.

Pour le registre injonctif, les régulations interactives peuvent être associées à une information relative au contexte de la séance tels que l'habillage de l'énoncé, les consignes sur le déroulement de la séance et les modalités temporelles, et à une action :

- d'explicitation ou de modification des critères d'évaluation et des buts de l'activité,

- de modification de l'habillage de l'énoncé,

- de validation,

- de relance ou reformulation de la question ou de la réponse de l'élève,

- de demande auprès d'un·e élève d'aider un·e autre élève ;

- de comparaison entre la production de deux élèves ou entre la production d'un·e et un référentiel.

Les actions définies ci-dessus sont des actions indépendantes de la tâche mathématique proposée et où le savoir enseigné ou prescrit est pris comme référentiel sans accompagnement de l'élève dans sa progression de résolution de la tâche mathématique.

Pour le registre affectif, les régulations interactives peuvent être associées à une information relative aux traits affectifs locaux de l'élève – eg. émotions, attitudes ou motivation –, aux traits affectifs globaux de l'élève – eg. estime de soi, croyances ou représentation – ou à la capacité de l'élève à réguler ses émotions (agir par curiosité, considérer la frustration comme un signal pour modifier sa stratégie, mobiliser des émotions positives pour s'engager dans la tâche mathématique, atténuer ou supprimer les émotions négatives pour se concentrer sur la tâche mathématique, etc.). À partir de l'interprétation de ces informations le ou la professeur·e peut agir sur les compétences affectives et les structures affectives des élèves (Debellis et Goldin, 2006), ainsi qu'agir sur l'estime de soi, l'auto-efficacité et les croyances des élèves vis-a-vis des mathématiques (Hattie et Timperley, 2007).

Chaque action menée par le ou la professeur·e peut s'effectuer avec des degrés de conscience et d'intentionnalité variables pouvant déboucher sur un changement ou un renforcement des mécanismes de régulation de l'élève (Allal, 2007). L'état des mécanismes de régulation de l'élève et de leurs évolutions par l'intervention du professeur·e sont importants à prendre en compte dans le cas de cet étude comparative. En effet, la comparaison des types de régulation des élèves permettra d'observer quels « schémas d'interaction » favorisent ou non la transformation des compétences de l'élève selon la tâche mathématique et s'il y a une variabilité selon le genre de l'élève. Les régulations interactives seront considérées comme efficaces si elles permettent :

> « le guidage, le contrôle ou l'ajustement des activités cognitives, affectives et sociales, favorisant ainsi la transformation des compétences de l'apprenant » (Allal, 2007).

Dans le cas contraire, elles seront considérées comme inefficaces. P. A. Genoud et M. Guillod (2014) ont montré à l'aide d'un questionnaire évaluant les attitudes socio-affectives en maths – 416 élèves y ont participés - que les croyances des élèves sont fortement liées à leurs émotions qui influencent l'activation des processus cognitifs. Ce constat montre la pertinence à porter un intérêt aux effets des « schémas d'interaction » sur l'élève et sur les compétences à acquérir, comme indicateur de la menace du stéréotype de genre sur l'apprentissage des mathématiques des élèves.

## V.    CONCLUSION

Cet article, rédigé dans le cadre de mes premiers travaux de thèse en didactique des mathématiques, discute du besoin d'appréhender ce qui se joue réellement au niveau de la performance des élèves lors d'une activité mathématique. De nombreuses études qualitatives sont menées sur la menace du stéréotype de genre sur l'apprentissage des mathématiques, l'influence du stéréotype du genre sur les pratiques évaluatives des professeur·es et le choix d'orientation des élèves selon leur genre. Néanmoins, très peu d'entre elles font une analyse par rapport à une compétence ou une tâche mathématique spécifique. Ce constat m'a amené à croiser différents outils théoriques de la didactique des mathématiques, de la psychologie cognitive et de l'affect mathématique pour m'inscrire dans le double champs de la recherche sur l'influence des stéréotypes de genre sur les pratiques évaluatives des professeur·es et l'effet de menace du stéréotype de genre sur la performance des élèves. Mon projet de recherche vise à apporter des éléments de réponse quant aux liens de dépendance entre les pratiques évaluatives des professeur·es, la performance des élèves et leur affect mathématique, afin d'évaluer l'existence d'un contrat didactique en évaluation différencié selon le sexe de l'élève.

REFERENCES

Allal, L. (2007). Introduction. Régulations des apprentissages : Orientations conceptuelles pour la recherche et la pratique en éducation. In *Régulation des apprentissages en situation scolaire et en formation* (p. 7-23). De Boeck Supérieur; Cairn.info.

Barrier, T., Hache, C., Mathé, A.-C., & Montigny, S. (2013). *Décrire l'activité géométrique des élèves : Instruments, regards, langage.* COPIRELEM, Nantes, France.

Bishop, A. (2014). Values in Mathematics Education. In S. Lerman (Éd.), *Encyclopedia of Mathematics Education* (p. 633-636). Springer Netherlands.

Croguennec, F. (2023). Résolution de situations-problèmes au primaire : Un défi de différenciation pédagogique et de didactique. *Revue des sciences de l'éducation*, *49*(1).

Darragh, L., & Radovic, D. (2018). Mathematics Learner Identity. In S. Lerman (Éd.), *Encyclopedia of Mathematics Education* (p. 1-4). Springer International Publishing.

Debellis, V. A., & Goldin, G. A. (2006). Affect and Meta-Affect in Mathematical Problem Solving : A Representational Perspective. *Educational Studies in Mathematics*, *63*(2), 131-147. JSTOR.

DEPP. (2024, juin). *Les évaluations nationales exhaustives* [Site gouvernemental]. Education gouv. https://www.education.gouv.fr/les-evaluations-nationales-exhaustives-307627

Fischer, J.-P., & Thierry, X. (2022). Boy's math performance, compared to girls', jumps at age 6 (in the ELFE's data at least). *British Journal of Developmental Psychology*, *40*, 504-519.

Genoud, P. A., & Guillod, M. (2014). Développement et validation d'un questionnaire évaluant les attitudes socio-affectives en maths. *Recherches en éducation*, *20*.

Grapin, N., & Mounier, E. (2020). Point de vue didactique sur les évaluations nationales françaises au début de la scolarité obligatoire. *Revue de Mathématiques pour l'École*, *234*, 22-30.

Grapin, N., & Mounier, É. (2024). Évaluer la résolution de problèmes additifs au CP : une étude exploratoire autour de la taille des nombres. *Grand N*, *113*, 29-48.

Hattie, J., & Timperley, H. (2007). The Power of Feedback. *Review of Educational Research*, *77*(1), 81-112.

Imberdis, A., Toczek, M.-C., Dutrévis, M., & Sacré, M. (2021). Disciplines scolaires et stéréotypes de genre : Perceptions d'élèves et d'enseignant·es. *L'orientation scolaire et professionnelle*, *50*(4).

Jamain, L. (2019). *Biais d'auto-évaluation de compétence en français et en mathématiques chez les élèves de primaire : Évolution et implications pour l'adaptation psychosociale et la réussite scolaire des élèves ?* [Éducation]. Université de Grenoble Alpes.

Kiwan-Zacka, M., & Roditi, É. (2018). Régulations didactiques et apprentissages des élèves. *Espace Mathématique Francophone 2018*, 77-85.

Martinot, D. (2001). Connaissance de soi et estime de soi : Ingrédients pour la réussite scolaire. *Revue des sciences de l'éducation*, *27*(3), 483-502. Érudit.

Mathieu-Benmerah, M. (2023). Épisodes évaluatifs et égalité. Comparaison didactique en mathématiques et EPS à l'école primaire. *Genre Éducation Formation (GEF)*, *7*.

Mella, N., & al. (2021). Socio-Emotional Competencies and School Performance in Adolescence : What Role for School Adjustment? *Frontiers in Psychology*, *12*(640661). Archive ouverte UNIGE.

Perronnet, C., Marc, C., & Paris-Romaskevich, O. (2024). *Matheuses. Les filles, avenir des mathématiques*. CNRS Éditions.

Petitfour, E. (2015, août). *Catégorisation des aides apportées dans une construction instrumentée*. EEDM18 : 18e école d'été de didactique des mathématiques, Brest, France.

Sayac, N. (2017). *Approche didactique de l'évaluation et de ses pratiques en mathématiques : Enjeux d'apprentissages et de formation* [Education]. Université Paris Diderot - Paris 7.

Sayac, N. (2024). Comment expliquer les écarts de performance entre les filles et les garçons en mathématiques aux évaluations nationales de CP-CE1 ? *Évaluer - Journal international de recherche en éducation et formation*, *9*(3), 32-51.

Veyrunes, P., Durny, A., Flavier, E., & Durand, M. (2005). L'articulation de l'activité de l'enseignant et des élèves pour résoudre un problème de mathématiques à l'école primaire : Une étude de cas. *Revue des sciences de l'éducation*, *31*(2), 471-489.

ANNEXES

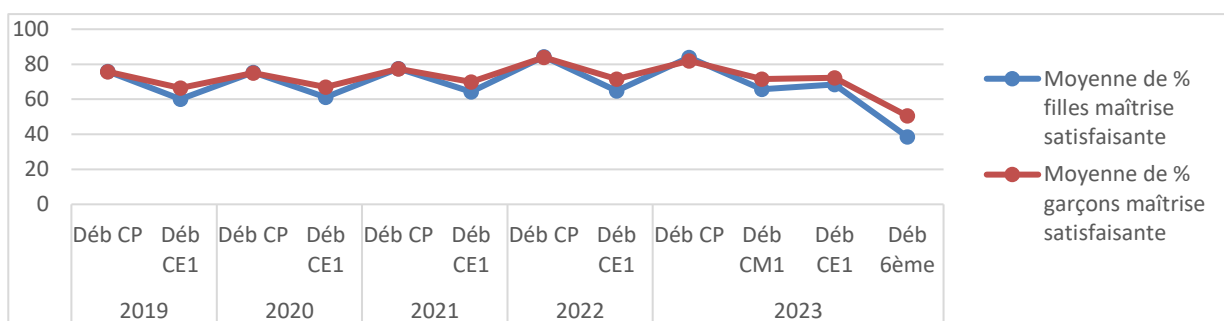

**Graphique 1** - *Écarts de performance filles–garçons en maîtrise satisfaisante en « Nombres et Calculs » de 2018 à 2024*

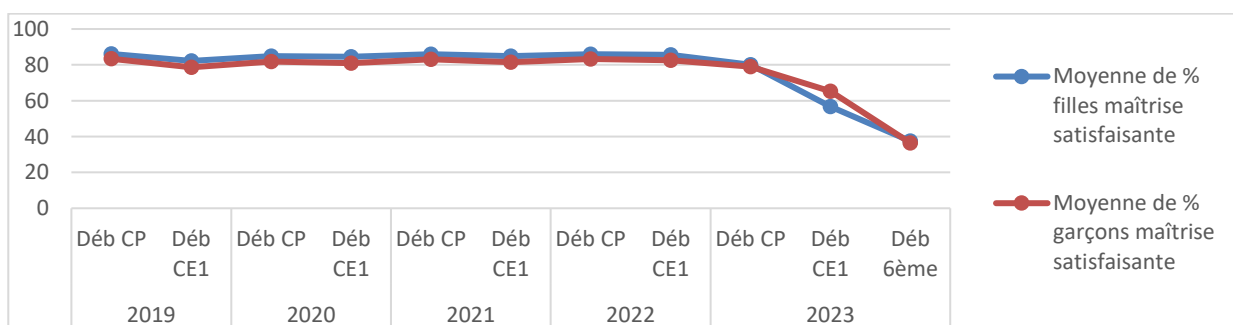

**Graphique 2** - *Écarts de performance filles–garçons en maîtrise satisfaisante en « Espace et Géométrie » de 2018 à 2023*

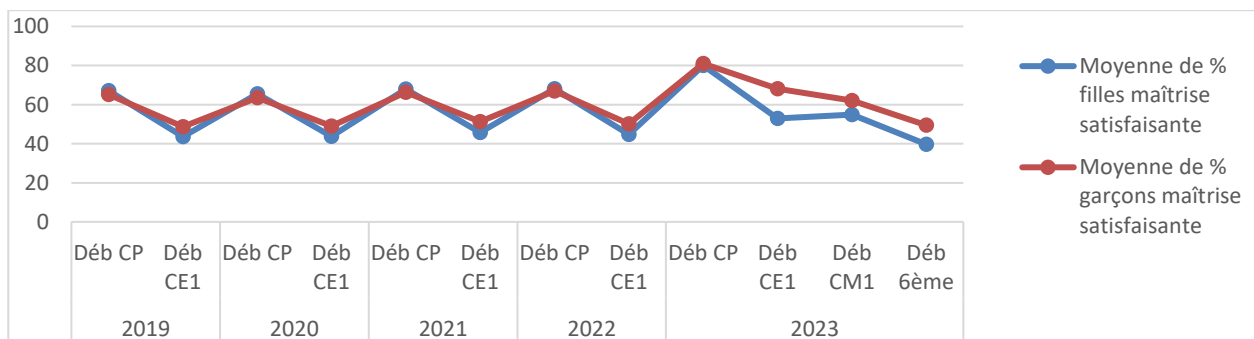

**Graphique 3** - *Écarts de performance filles–garçons en maîtrise satisfaisante en « Résolution de problèmes » de 2019 à 2024*